\documentclass[12pt,reqno,intlimits]{amsart}
\usepackage{amsfonts,amsmath,amsxtra,amsthm,amssymb,latexsym,pb-diagram}

\usepackage[cp1251]{inputenc}
\usepackage[T2A]{fontenc}
\usepackage{graphicx, epsfig}
\usepackage{mathrsfs}
\theoremstyle{plain}
\newtheorem{theorem}{Theorem}[section]
\newtheorem{corollary}[theorem]{Corollary}
\newtheorem{proposition}[theorem]{Proposition}

\theoremstyle{definition}
\newtheorem{definition}[theorem]{Definition}

\begin{document}

\title{Sublinear extension of Grothendieck Algebraic $K$--theory}
{\it{This paper is dedicated to the 80-th anniversary of Yuri I. Manin}}

\author{Igor V. Orlov}
\maketitle



\smallskip
\begin{abstract}{A commutative diagram that connects the basic objects of commutative algebra with the main objects of commutative analysis is constructed. Namely, with the help of five types of canonical embeddings we constructed a diagram between two sets of objects: Abelian semigroups -- Abelian regular (cancellative) semigroups -- Abelian groups, on the first hand, and convex cones -- regular convex cones -- linear spaces, on the other hand. Thus, some extension of the Grothendieck algebraic $K$--theory arises, that includes the basic objects not only of linear (smooth) analysis but of sublinear (nonsmooth) analysis also.}\end{abstract}

\section{Basic objects: terminology, notation, auxiliary results}

\begin{enumerate}
  \item \emph{Abelian semigroup } with zero element (i.e., monoid): $X=\{x\},$ with additive notation. The corresponding category denote by $(S)$.

  \item \emph{Regular (cancellative) Abelian semigroup} is an Abelian semigroup that satisfies the \emph{ cancellation law}: $(x+y=x+z)\Rightarrow(y=z)$ (see \cite{Clifford}). The corresponding category denote by $(RS)$.
  \item \emph{Abelian group.}  The corresponding category denote by $(G)$.

  \item \emph{Convex cone} is an Abelian semigroup with respect to vector addition that forms a module over $\mathbb{R}_+$ with respect to multiplication by scalars (see \cite{Polovinkin, Orlov}). The corresponding category denote by $(Con)$.

  \item \emph{Regular (cancellative) convex cone} is a convex cone that satisfies the cancellation law (see \cite{Polovinkin, Orlov}). The corresponding category denote by $(RCon)$.
  \item \emph{Linear space} (over $\mathbb{R}$). The corresponding category denote by $(Lin)$.
\end{enumerate}
\vspace{6pt}

Let's remind some auxiliary concepts and results.

  \begin{proposition}{\rm{(\cite{Clifford})}} An arbitrary Abelian semigroup can be isomorphically (injectively and additively) embedded into some Abelian group if and only if it is regular.
  \end{proposition}

  \begin{corollary}{\rm{(\cite{Polovinkin})}} An arbitrary convex cone can be isomorphically (i.e. injectively and $\mathbb{R}_+-$linearly) embedded into some linear space if and only if it is regular.
  \end{corollary}

   \begin{definition}{(\cite{Manin,Atiyah,Achar})} The minimal Abelian group that contains given regular semigroup $X$ is called \emph{Grothendieck group} of $X$ and is denoted by $Gr(X).$ Respectively, the minimal linear space that contains given regular convex cone $X$ let's call \emph{Grothendieck linear group} of $X$ and is denoted by $Gr_L(X).$
   \end{definition}

 \section{Divisible objects: terminology, notation, auxiliary results}

\begin{enumerate}
  \item \emph{Divisible Abelian semigroup} is an Abelian semigroup $X$ that satisfies the condition:
      $$\forall\, x\in X\,\, \forall\, n\in \mathbb{N}\,\,\exists\, y\in X:  \underbrace{y+ \dots + y}_n:={\sum}_{n}y=x$$
(see \cite{Griffith}). The corresponding categories denote by $(DS)$,  $(DRS)$, $(DG)$.

\medskip
  \item {Uniquely divisible Abelian semigroup} is a divisible Abelian semigroup $X$ that satisfies the condition:
      $$\forall\, x_1, x_2\in X,\ \forall\, n \in \mathbb{N} \,\, \left({\sum}_{n}x_1={\sum}_{n}x_2\right)\Rightarrow (x_1=x_2)$$
      (see \cite{Griffith, Tabor}). The corresponding categories denote by $(US)$,  $(URS),\ (UG)$.
\end{enumerate}

\begin{proposition} Let $X$ be a uniquely divisible Abelian semigroup. Then:
  $$\forall\, x\in X,\,\,\, \forall\, n_1,n_2\in \mathbb{N}\quad \left({\sum}_{n_1}x={\sum}_{n_2}x, \,x\neq0\right)\Rightarrow (n_1=n_2).$$
  \end{proposition}

  \section{Basic canonical embeddings}

  Let's describe in short the main embeddings that will be used further.

\smallskip
\begin{enumerate}

 \item\emph{Regularization.} $R:\,(S)\to(RS)$, $(DS)\to (DRS)$, $(US)\to(URS)$, $(Con)\to (RCon).$

\smallskip
 Factorize  $X$ by $ (y_1Ry_2)\Leftrightarrow(\exists\, x\in X:\, x+y_1=x+y_2).$

Let  $X_R=X/ R$ be corresponding factor semigroup, then $R:\,X\to X_R$, $y\mapsto\{z\in X|\ yRz\}$ is the required embedding and $X_R\in (RS)$.

\smallskip
\item\emph{Formal difference.} $F:\,(RS)\to(G)$, $(DRS)\to (DG)$, $(URS)\to(UG)$, $(RCon)\to (Lin).$

\smallskip
Factorize $X\times X$ by \[\left((y_1,z_1)F(y_2,z_2)\right) \Leftrightarrow \left( y_1+z_2=y_2+z_1\right) .\]
Let $X_F=X\times X/ F$ be corresponding factor semigroup. Introduce the subtraction operation in $X_F$ by involution
\[-(y,z)=(z,y).\] Then  $F:\,X\to X_F$ $\left(x\mapsto \{(y,z)\big|\ x+y=z\}\right)$ is the  required embedding and   $X_F=F(X)-F(X).$

\smallskip
\item\emph{Divisibility.} $D:\,(S)\to(DS)$, $(RS)\to (DRG)$, $(G)\to(DG).$

\smallskip
Factorize $X\times \mathbb{N}$ by
$$\left((x_1,n_1)D(x_2,n_2)\right) \Leftrightarrow \left( {\sum}_{n_2}x_1={\sum}_{n_1}x_2\right). $$
Let $X_D=X\times \mathbb{N}/ D$ be corresponding factor semigroup, then $D:\,X\to X_D$
$\left(x\mapsto \{(y,n)\big|\ (x,1)D(y,n)\}\right)$ is the required embedding.

\smallskip
\item\emph{Uniquely divisibility.} $U:\,(DS)\to(US)$, $(DRS)\to (URG)$, $(DG)\to(UG).$

Factorize $X=$ by $ (x_1Ux_2)\Leftrightarrow\left(\exists\, n\in \mathbb{N}: \ {\sum}_{n}x_1={\sum}_{n}x_2\right).$

Let $X_U=X/ U$ be corresponding factor semigroup, then $U:\,X\to X_U$
$(x\mapsto \{y\in X|\,\exists\, n\in \mathbb{N}:\,{\sum}_{n}x= {\sum}_{n}y\})$ is the required embedding.

\smallskip
\item\emph{Modulation.} (see \cite{Orlov}) $M:\,(US)\to(Con)$, $(URS)\to (RCon)$, $(UG)\to(Lin).$

\smallskip
Let's introduce in uniquely divisible Abelian semigroup $X$ an ``additive multiplication'' \ by non-negative scalars, first for rational case.

\begin{enumerate}
  \item For $x\in X,\ r=\frac{m}{n}\in \mathbb{Q}_+$ set $(y=r*x)\Leftrightarrow(\sum_{m}x=\sum_{n}y).$
  \item For $\gamma\in \mathbb{R}_+$ that defined by Dedekind cutting $A|B$ in $ \mathbb{Q}_+,$ set \[\gamma*x=(A*x\big|B*x).\]
  \item Define $X_M$ as additive envelope of the set  $\mathbb{R}_+* X$ with respect to Minkowsky addition:
      \[X_M=\left\{\sum\limits_{k=1}^n\gamma_k*x_k\big|\,\gamma_k\in \mathbb{R}_+,\ x_k\in X,\ n\in \mathbb{N}\right\}.\]  Here the modulation is extended to $X$  by obvious way:
      \[\alpha*(\sum_{k=1}^n\gamma_k*x_k)=\sum_{k=1}^n(\alpha\gamma_k)*x_k.\]
\end{enumerate}

The canonical embedding $M:\,X\to X_M$ is defined by the equality:
\[Mx=1*x=\left([0,1]_\mathbb{Q}*x\big|(1,+\infty)_\mathbb{Q}*x\right).\]

\end{enumerate}

  \section{Properties of the basic embeddings}

  The following statements can be checked by direct transformations.

      \begin{theorem} The embedding $R$ is an additive homomorphism from $X$ \emph{onto} $X_R$.
      \end{theorem}

 \begin{theorem} The embedding $F$ is an additive isomorphism from $X$ \emph{into} $X_F$. In addition, in case  $(RCon)\to (Lin),\ F$ is $\mathbb{R}_+$-linear isomorphism.
 \end{theorem}

 \begin{theorem} The embedding $D$ is an additive isomorphism from $X$ \emph{into} $X_D$. In addition,  $X_D=D(X)-D(X).$
 \end{theorem}

\begin{theorem}  The embedding $U$ is an additive homomorphism from $X$ \emph{onto} $X_U$.
\end{theorem}

\begin{theorem} The embedding $M$ is an additive isomorphism from $X$ \emph{into} $X_M$. In addition, $X_M=Add(\mathbb{R}_+*X).$
\end{theorem}

\newpage
  \section{The main result}
      \begin{theorem} The following diagram is commutative, together with any its subdiagram.
      \end{theorem}

 \begin{figure}[h]
\center \includegraphics[width=300pt]{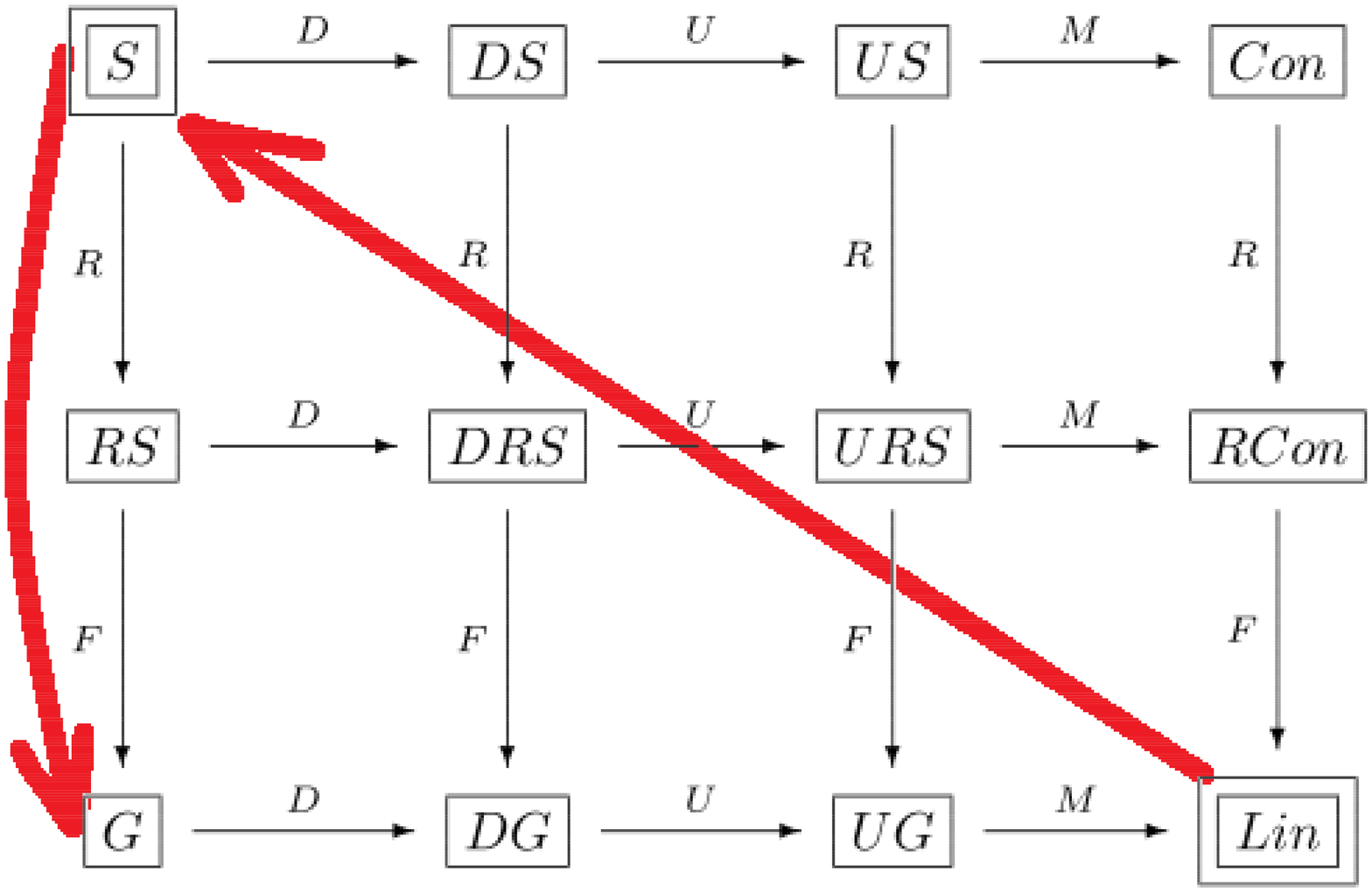}
\label{2p}
\end{figure}

\noindent{\bf{Final Remark.}} The red arrows on the diagram below show a connection between Theorem~5.1 and classical Grothendieck Theorem (see \cite{Manin},\cite{Atiyah},\cite{Achar},\cite{Karoubi},\cite{Weibel}).

\bigskip
Igor Vladimirovich Orlov

Department of Mathematics and Informatics

Crimean Federal V. Vernadsky University

Academician Vernadsky Ave., 4

Simferopol, Republic Crimea, Russia, 295007

E-mail:   igor\_v\_orlov@mail.ru\\
\end{document}